\documentclass{amsart}
\usepackage{amscd}
\usepackage{amssymb}
\usepackage[curve]{xypic}
\newbox\mybox
\def\overtag#1#2#3{\setbox\mybox\hbox{$#1$}\hbox to
  0pt{\vbox to 0pt{\vglue-#3\vglue-\ht\mybox\hbox to \wd\mybox
      {\hss$\ss#2$\hss}\vss}\hss}\box\mybox}
\def\undertag#1#2#3{\setbox\mybox\hbox{$#1$}\hbox to 0pt{\vbox to
    0pt{\vglue#3\vglue\ht\mybox\hbox to \wd\mybox
      {\hss$\ss#2$\hss}\vss}\hss}\box\mybox}
\def\lefttag#1#2#3{\hbox to 0pt{\vbox to 0pt{\vss\hbox to
      0pt{\hss$\ss#2$\hskip#3}\vss}}#1}
\def\righttag#1#2#3{\hbox to 0pt{\vbox to 0pt{\vss\hbox to
      0pt{\hskip#3$\ss#2$\hss}\vss}}#1}
\let\ss\scriptstyle

\def\Dot{\lower.2pc\hbox to 2pt{\hss$\bullet$\hss}}
\def\Circ{\lower.2pc\hbox to 2pt{\hss$\circ$\hss}}
\def\Vdots{\raise5pt\hbox{$\vdots$}}
\newcommand\lineto{\ar@{-}}
\newcommand\dashto{\ar@{--}}
\newcommand\dotto{\ar@{.}}

\newcommand\valence{valency}
\newcommand\Brieskorn{Brieskorn-Pham}

\newcommand{\comment}[1]{\marginpar{\sffamily{\tiny #1
\par}\normalfont}}
\renewcommand{\comment}[1]{}
\newcommand\quotient[2]{#1/#2}

\newcommand\SU{\operatorname{SU}}
\newcommand\SO{\operatorname{SO}}

\newcommand\Q{{\mathbb Q}}
\newcommand\R{{\mathbb R}}
\newcommand\C{{\mathbb C}}
\newcommand\Z{{\mathbb Z}}
\newcommand\N{{\mathbb N}}
\renewcommand\O{{\mathcal O}}
\newtheorem{theorem}{Theorem}[section]
\newtheorem*{theorem*}{Theorem}

\newtheorem{proposition}[theorem]{Proposition}

\newtheorem{conjecture}{Conjecture}

\newtheorem*{conjecture*}{Conjecture}
\theoremstyle{definition}
\newtheorem{definition}[theorem]{Definition}
\newtheorem{example}{Example}
\newtheorem*{example1}{Example 1 \rm(ctd.)}
\newtheorem{examples}[example]{Examples}

\newtheorem*{example*}{Example}
\newtheorem*{remark*}{Remark}

\begin{document}
\title{Universal abelian covers of surface singularities}
\author{Walter D. Neumann}
\thanks{Research supported under NSF grant no.\ DMS-0083097}
\address{Department of Mathematics\\Barnard College, Columbia
  University\\New York, NY 10027}
\email{neumann@math.columbia.edu}
\author{Jonathan Wahl}
\thanks{Research supported under NSA grant no.\ MDA904-98-1-0016}
\address{Department of Mathematics\\The University of
North Carolina\\Chapel Hill, NC 27599-3250}
\email{jw@math.unc.edu}
\begin{abstract}
  We discuss the evidence for and implications of a conjecture that
  the universal abelian cover of a $\Q$-Gorenstein surface singularity
  with finite local homology (i.e., the singularity link is a
  $\Q$-homology sphere) is a complete intersection singularity.
\end{abstract}
\maketitle

\section{History}

Our story starts with a result of Felix Klein.  The \emph{icosahedral
  group} $I\subset \SO(3)$ is the group of orientation preserving
symmetries of a regular icosahedron in $\R^3$. The binary icosahedral group
$I'\subset \SU(2)$ is the inverse image of $I$ under the double cover
$\SU(2)\to\SO(3)$. The quotient $\quotient{\C^2}{I'}$ is a complex
variety with an isolated singularity at $0$.
Felix Klein showed that
$$\quotient{\C^2}{I'}\cong V(2,3,5),$$
where $V(p,q,r)$ is the so-called \emph{\Brieskorn{} variety}
$$V(p,q,r):=\{(z_1,z_2,z_3)\in\C^3 :
z_1^p+z_2^q+z_3^r=0\}.$$

As the reader probably knows, the binary icosahedral group is
intimately related to the triple $(2,3,5)$. The non-cyclic finite
subgroups of $\SU(2)$ are the binary dihedral and polyhedral groups,
and they are classified by the \emph{Platonic triples} $(p,q,r)$:
triples of integers satisfying
$$1<p\le q\le r,\quad\frac1p+\frac1q+\frac1r>1.$$
These are the triples
$$(2,2,n),~n\ge2;\quad(2,3,3),(2,3,4),(2,3,5).$$
For a Platonic triple
$(p,q,r)$ the group of isometries of the 2-sphere generated by
reflections in the sides of a spherical triangle with angles $\frac\pi
p,\frac\pi q, \frac\pi r$ is the \emph{$(p,q,r)$ triangle group}.  Its
orientation preserving subgroup $\Delta(p,q,r)$ is the dihedral group
of order $2n$ for $(p,q,r)=(2,2,n)$, while $\Delta(2,3,3)$,
$\Delta(2,3,4)$ and $\Delta(2,3,5)$ are the tetrahedral, octahedral
and icosahedral groups. They are subgroups of $\SO(3)$, and on lifting
to $\SU(2)$ we get the corresponding binary groups $\Delta'(p,q,r)$;
these give all finite non-cyclic subgroups of $\SU(2)$.

It is thus natural to guess that the relationship
$V(2,3,5)\cong\quotient{\C^2}{\Delta'(2,3,5)}$ extends to the other
Platonic triples, but this is not correct, as the following table
shows:
\begin{align*}
V(2,3,5)&\cong\quotient{\C^2}{\Delta'(2,3,5)}\\
V(2,3,4)&\cong\quotient{\C^2}{\Delta'(2,3,3)}\\
V(2,3,3)&\cong\quotient{\C^2}{\Delta'(2,2,2)}\\
V(2,2,n)&\cong\quotient{\C^2}{C_n}.
\end{align*}
(Here, $C_n$ is the cyclic group of order $n$.)

The correct relationship is, in fact,
$$V(p,q,r)\cong\quotient{\C^2}{[\Delta'(p,q,r),\Delta'(p,q,r)]},$$
where $[G,G]$ denotes the commutator subgroup of a group $G$. In other
words, $V(p,q,r)$ is the \emph{universal abelian cover} of
$\quotient{\C^2}{\Delta'(p,q,r)}$ (this is literally true after
one removes the singular
points from each space, where the map is ramified).

This relationship was shown to extend to the euclidean and hyperbolic
triangle groups by J. Milnor \cite{milnor} and then to groups arising
from general polygonal reflection groups in \cite{neumann77}. However,
the connection with reflection groups turned out to be a red herring,
as was  shown in \cite{neumann83}. The essential fact appeared
to be that the varieties in question were weighted homogeneous (i.e.,
admit effective $\C^*$-actions):
\begin{theorem}[\cite{neumann83}]\label{th:weighted homog}
  Let $(X,o)$ be a weighted homogeneous normal surface singularity
  whose link $M$ (boundary of a closed regular neighborhood) is a
  $\Q$-homology sphere (i.e., $H_1(M;\Q)=0$).  Then the universal
  abelian cover $(\tilde X,o)$ of $(X,o)$ is a \Brieskorn{} complete
  intersection.
\end{theorem}
\noindent
A \emph{\Brieskorn{} complete intersection} is a surface given by $n-2$
equations $$
a_{1j}z_1^{p_1}+\dots+a_{nj}z_n^{p_n}=0,\quad
j=1,\dots,n-2,$$
where $p_1,\dots,p_n$ are fixed integers $\ge2$ and
the $a_{ij}$ are sufficiently general coefficients.  The exponents
$p_i$ in this theorem are, in fact, the orders of the isotropy groups
of the action of $\C^*$ on $X$.

\medskip The story does not stop here.  One of us (J.W.) suspected that
``weighted homogeneous'' was also a red herring, and finally
persuaded the other that this might be so.  We were nevertheless
surprised at how general the phenomenon appears to be.

Before describing this in more detail, we would like to stress why
a result of this type should be surprising.  First, while $\tilde X$
in Theorem 1.1 is easily shown to be Gorenstein, it is well-known that being a
complete intersection is a much stronger property, and there are no general
results we know of which distingush the two. Second, among the
weighted homogeneous complete intersections, there are many which
are not even of the same topological type as a \Brieskorn{} one.

There are two obvious necessary conditions for the universal
abelian cover of a singularity to be a complete intersection. In
this situation, we need to consider $(X,o)$ as a germ of a
singularity; we may thus assume $X$ is a regular neighborhood of
$o$. The maximal unramified abelian cover of $X-\{o\}$ is a Galois
cover with covering transformation group $H_1(M;\Z)$, where $M$ is
the link of the singularity. If this is a finite cover, then it can
always be completed by adding a single point to give a new germ and a
map
$(\tilde X,o)\rightarrow (X,o)$,
which is what we mean by the \emph{universal
  abelian cover of $(X,o)$}. The two necessary conditions are:
\begin{enumerate}
\item \label{it:exist} ($\Q$-homology sphere link). The universal abelian
  cover must exist.  That is, $H_1(M;\Z)$ must be finite, or $M$ is a
  $\Q$-homology sphere.
\item \label{it:q-gorenstein} ($\Q$-Gorenstein). $(X,o)$ should be
  \emph{$\Q$-Gorenstein}, that is the dualizing sheaf $\omega_X$ has
  finite order, or, equivalently, for some $n>0$, $\Gamma(X-\{o\},
  \omega_X^{\otimes n})$ is a free $\O_X$-module.  Alternatively,
  $(X,o)$ is the quotient of some Gorenstein singularity by a finite
  group acting freely off the singular point.
\end{enumerate}
 It should also be noted that the Gorenstein condition not only
limits the topology, but in the ``moduli space'' of all
singularities with a given topology, non-Gorenstein singularities
are nearly always generic (Laufer, in \cite{laufermin}, shows the only
exceptions are rational double points or minimally elliptic
singularities).

\begin{conjecture}\label{conj:main}
  Let $(X,o)$ be a $\Q$-Gorenstein singularity whose link is a
  $\Q$-homology sphere.  Then the universal abelian cover of $(X,o)$
  is a complete intersection singularity.
\end{conjecture}

We have already said why this conjecture is surprising. In the
next section we describe what the complete intersection equations
conjecturally look like; these ``splice diagram equations''
generalize the \Brieskorn{} part of Theorem 1.1. The possible
topologies of normal surface singularities have been known for
decades (through the work of Grauert and Waldhausen, see e.g.,
\cite{neumann81}), but it is nevertheless an extremely difficult
problem in general to give an explicit algebraic description for a
singularity with given topology. It has so far been done only for
limited types of topology.  In contrast, our conjectures would
give explicit algebraic descriptions for a remarkably broad class
of topologies.

The next section gives a series of (increasingly speculative)
conjectures that refine our main Conjecture \ref{conj:main}. Fuller
details concerning these conjectures and partial results described
here will be published elsewhere.\comment{WDN changed this paragraph}

\section{The Conjectures}\label{sec:conjectures}

Let $(X,o)$ be a normal surface singularity germ and $M$ its link,
that is, the boundary of a regular neighborhood of $o$ in $X$. We
shall assume throughout this section that $M$ is a $\Q$-homology
sphere.

Let $\pi\colon\overline X\to X$ be a good resolution. ``Good''
means that the exceptional divisor $D=\pi^{-1}(o)$ has only normal
crossings.  The $\Q$-homology sphere condition is equivalent to
$D$ being \emph{rationally contractible}; that is,
\begin{itemize}
\item each component of $D$ is a smooth rational curve;
\item the dual resolution graph $\Gamma$ (the graph with a vertex for
  each component of $D$ and an edge for each intersection of two
  components) is a tree.
\end{itemize}

We weight each vertex $v$ of $\Gamma$ by the self-intersection
number $E_v\cdot E_v$ of the corresponding component $E_v$ of $D$. The
\emph{intersection matrix} for $\Gamma$ is the matrix $A(\Gamma)$
with entries $a_{vw}=E_v\cdot E_w$, that is,
\begin{align*}
  a_{vw}&=1\quad\text{if $v\ne w$ and $v$ and $w$ are joined by an
    edge}\\
  a_{vw}&=0\quad\text{if $v\ne w$ and $v$ and $w$ are not joined by an
    edge}\\
a_{vv}&=E_v\cdot E_v
\end{align*}
It is well known that $A(\Gamma)$ is negative-definite and its
cokernel (also called the \emph{discriminant
group}) is $H_1(M)$.
In particular, $$d(\Gamma):=\det(-A(\Gamma))$$
is the order of $H_1(M)$.
By a
fundamental result of \cite{neumann83}, with a few exceptions
$A(\Gamma)$ determines and is determined by the topology of M.

To write down equations of the universal abelian cover, we associate
another combinatorial object to the weighted dual graph $\Gamma$---a
\emph{splice diagram} \cite{eisenbud}.  When $A(\Gamma)$ is unimodular
(i.e., M is an \emph{integral} homology sphere), this association can
be found in \cite{eisenbud}. In that situation one can also recover
the resolution graph $\Gamma$ from the splice diagram, but without
unimodularity one usually cannot --- many resolution graphs can share
the same splice diagram.

\begin{definition}[Splice diagram of a rationally contractible graph]
  Let $\Delta$ be the tree obtained by replacing each maximal string
  in $\Gamma$ by a single edge (a \emph{string} is a simple path in
  $\Gamma$ whose interior is open in $\Gamma$). Thus $\Delta$ is
  homeomorphic to $\Gamma$ but has no vertices of \valence{} $2$. At each
  \emph{node} $v$ (vertex of \valence{} $\ge 3$) of $\Delta$ we introduce
  weights on the incident edges $e$ as follows. Let $\Gamma_{ve}$ be the
  subgraph of $\Gamma$ cut off by the edge of $\Gamma$ at $v$
  in the direction of $e$ as in the following picture. The corresponding
  weight is then $d_{ve}:=d(\Gamma_{ve})$.
    $$
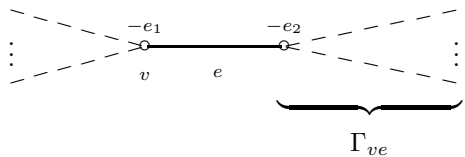
\xymatrix@R=6pt@C=24pt@M=0pt@W=0pt@H=0pt{
&&&&&&\\
\vdots&&\undertag{\overtag{\Circ}{-e_1}{8pt}}{v}{8pt}\lineto[rr]\dashto[ull]\dashto[dll]&{.}\undertag{}{e}{8pt}&
  \overtag{\Circ}{-e_2}{8pt}\dashto[urr]\dashto[drr]&&\vdots\\
&&&&&&\\
&&&&&{\hbox to 0pt{\hss$\underbrace{\hbox to 70pt{}}$\hss}}&\\
&&&&&{\Gamma_{ve}}&}$$
\end{definition}

For example, the following is a resolution graph and its splice diagram.
$$
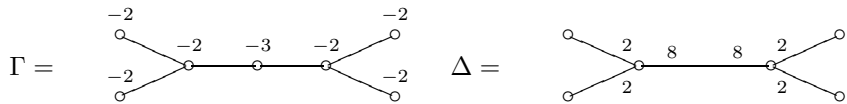
\xymatrix@R=6pt@C=24pt@M=0pt@W=0pt@H=0pt{
\\
&\overtag{\Circ}{-2}{8pt}&&&&\overtag{\Circ}{-2}{8pt}\\
{\Gamma=}&&\overtag{\Circ}{-2}{8pt}\lineto[ul]\lineto[dl]\lineto[r]&
\overtag{\Circ}{-3}{8pt}&\overtag{\Circ}{-2}{8pt}\lineto[ur]\lineto[dr]\lineto[l]&\\
&\overtag{\Circ}{-2}{8pt}&&&&\overtag{\Circ}{-2}{8pt}}
\qquad
\xymatrix@R=6pt@C=24pt@M=0pt@W=0pt@H=0pt{
\\
&\Circ&&&&\Circ\\
{\Delta=}&&\Circ\lineto[ul]_(.25){2}\lineto[dl]^(.25){2}\lineto[rr]^(.25){8}^(.75){8}&
&\Circ\lineto[ur]^(.25){2}\lineto[dr]_(.25){2}&\\
&{\Circ}&&&&{\Circ}\\&~}
$$
The weight $8$ on the left node of $\Delta$ is $d(\Gamma_{ve})$ with
$$
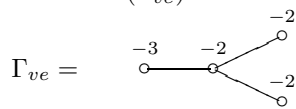
\xymatrix@R=6pt@C=24pt@M=0pt@W=0pt@H=0pt{
&&&\overtag{\Circ}{-2}{8pt}\\
{\Gamma_{ve}}=&\overtag{\Circ}{-3}{8pt}\lineto[r]&
\overtag{\Circ}{-2}{8pt}\lineto[ur]\lineto[dr]\\
&&&\overtag{\Circ}{-2}{8pt}}
$$

We will associate equations to every node of the splice diagram; some
preparation is needed.  An \emph{end} of $\Delta$ is a vertex of
\valence{} 1.  We first associate a weight $\ell_{vw}\in\N$ to a pair
consisting of a node $v$ and end $w$.  The weight $\ell_{vw}$ is the
product of the edge weights in $\Delta$ that are adjacent to, but not
on, the path from $v$ to $w$ in $\Delta$. Let $\ell'_{vw}$ be the same
product but excluding weights adjacent to $v$ itself. Thus,
$$\ell_{vw}d_{ve}=\ell'_{vw}d_v,$$
where $d_v$ is the product of edge
weights adjacent to $v$ and $e$ is the first edge of the path from $v$
to $w$.

\begin{example}\label{ex:cont}
For example, in the above example with vertices labelled as follows
$$
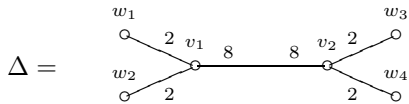
\xymatrix@R=6pt@C=24pt@M=0pt@W=0pt@H=0pt{
\\
&\overtag{\Circ}{w_1}{8pt}&&&&\overtag{\Circ}{w_3}{8pt}\\
{\Delta=}&&\overtag{\Circ}{v_1}{8pt}\lineto[ul]_(.45){2}\lineto[dl]^(.45){2}
\lineto[rr]^(.25){8}^(.75){8}&
&\overtag{\Circ}{v_2}{8pt}\lineto[ur]^(.45){2}\lineto[dr]_(.45){2}&\\
&\overtag{\Circ}{w_2}{8pt}&&&&\overtag{\Circ}{w_4}{8pt}\\&~}
$$
we have $\ell_{v_1w_1}=16$, $\ell'_{v_1w_1}=1$, $\ell_{v_1w_3}=8$,
$\ell'_{v_1w_3}=2$, etc.
\end{example}
We are interested in a numerical condition on the splice diagram.

\begin{conjecture}[Semigroup condition]\label{conj:sg}
  Let $\Delta$ be the splice diagram associated to a $\Q$-Gorenstein
  singularity with $\Q$-homology sphere link.  Then
  for each node $v$ and adjacent edge $e$ the edge-weight $d_{ve}$
  is in the semigroup
  $$\N\langle\ell'_{vw}:w \text{ an end of $\Delta$ in }
  \Delta_{ve}\rangle,$$
  where $\Delta_{ve}$ is the subtree of $\Delta$
  cut off by $e$.

  Equivalently, if $d_v$ is the product of the edge weights
  $d_{ve}$ adjacent to $v$, then
  $$d_{v}\in\N\langle\ell_{vw}:w \text{ an end of $\Delta$ in }
  \Delta_{ve}\rangle$$
\end{conjecture}
The condition of this conjecture is \emph{not} true for every
resolution tree.  In particular, this Conjecture would imply a
highly non-trivial (and unexpected) topological condition for a
singularity to be $\Q$-Gorenstein.  It is shown in
\cite{neumann-wahl01} that this condition is implied if certain
naturally defined knots in the link $M$ are cut out by hyperplane
sections of $(X,o)$.  The content of Conjecture 2 is that this
condition is necessary for the tree to be realizable by a
$\Q$-Gorenstein singularity.

So, let us now consider a splice diagram $\Delta$ satisfying the
semigroup condition (such as our example above).  The associated
complete intersection singularity will have embedding dimension less
than or equal to the number $n$ of ends of $\Delta$. Thus, to each end
$w$ we associate a variable $z_w$.

To each node $v$ of $\Delta$ of \valence{} $\delta_v$ we will
associate $\delta_v-2$ equations. It is easy to check that $n-2 =
\sum(\delta_v-2)$ (summed over the nodes of $\Delta$), so this
will give the right number of equations.

 We describe the $\delta_v-2$ equations
associated to a node $v$. For each end $w$ we give the variable
$z_w$ weight $\ell_{vw}$ (we call this the \emph{$v$-weight} of
$z_w$). To each edge $e$ at $v$ we can now associate one or more
\emph{admissible monomials} $\prod_w z_w^{\alpha_{vw}}$, product
over ends $w$ in $\Delta_{ve}$, with the exponents
$\alpha_{vw}\in\N$ satisfying
\begin{equation}\label{eq:sg}
    d_{v}=\sum_w\alpha_{vw}\ell_{vw},\quad\text{sum over $w$ an end in
      } \Delta_{ve}.
\end{equation}
Thus each admissible monomial has total $v$-weight $d_{v}$.  Each of
the $\delta_v-2$ equations associated to $v$ is obtained by equating
to zero some $\C$-linear combination of the admissible monomials
associated to $v$. If $\Delta$ has just one node this gives equations
of \Brieskorn{} type.

We choose $\delta_v-2$ equations in this way for each node $v$, giving
a total of $n-2$ equations in the $n$ variables $z_w$. We must also
require that the coefficients of the equations are ``sufficiently
general'' in an appropriate sense.  
\begin{definition}
  We say such a system of equations
is in \emph{splice diagram form}.
\end{definition}

\comment{Around here is rewritten. I've included stuff on suff.\ 
  general since it is really relevant to the example.}  We will not
specify in detail what ``sufficiently general'' means here. It is a
generic condition which usually depends on all the coefficients at
once. But if, for some node $v$, we have chosen just one admissible
monomial for each edge at $v$, and the $\delta_v-2$ equations
associated to $v$ involve only these $\delta_v$ monomials, then
``sufficiently general'' for these equations means that all maximal
minors of the $(\delta_v-2)\times\delta_v$ matrix of coefficients are
non-singular.  If $\Delta$ has one node (\Brieskorn{} complete
intersection) this is the well-known condition for an isolated
singularity due to H. Hamm.

Thus, for splice diagrams with one node, ``splice diagram form'' gives
exactly the isolated \Brieskorn{} complete intersection singularities.

\begin{example1}
  For the $\Delta$ of Example \ref{ex:cont} we abbreviate
  $z_{w_i}=z_i$.  The admissible monomials for vertex $v_1$ are
  $z_1^2$, $z_2^2$, and any $z_3^\alpha z_4^{\alpha'}$ with
  $2\alpha+2\alpha'=8$. Thus the equation associated to the node $v_1$
  might be $z_1^2+z_2^2+z_3^\alpha z_4^{4-\alpha}=0$, for some
  $0\le\alpha\le4$.  Similarly for the second node, so our system of
  equations might be
\begin{equation}\label{equations}
\begin{array}{r}
  z_1^2+z_2^2+z_3^\alpha z_4^{4-\alpha}=0,\\
  z_3^2+z_4^2+z_1^{\beta} z_2^{4-\beta}=0.
\end{array}
\end{equation}
This system is always in ``splice diagram form'' by our comments
above.  However, if we replace the equations, for example, by $az_1^2+bz_2^2+z_3^\alpha z_4^{4-\alpha}=0$ and
$z_3^2+z_4^2+cz_1z_2^3+dz_1^3z_2=0$, then the condition that the
coefficients are sufficiently general is: $ab\ne0$ and $ad-bc\ne0$.
\end{example1}

 \begin{theorem}
   A system of equations in splice diagram form defines an isolated
   complete intersection surface singularity.
\end{theorem}

We would like a complete intersection of this form to give the universal
abelian cover of our original singularity (up to equisingularity,
actually).  Since the covering
transformation group will be the discriminant group of $A(\Gamma)$,
this group should act in a natural way on the n-dimensional space
of the variables $z_w$.  This is indeed the case:

\begin{proposition}
  The discriminant group of $A(\Gamma)$ acts faithfully and diagonally
  on the n-dimensional space of the variables $z_w$, and acts freely
  in codimension $1$ (that is, no coordinate hyperplane is fixed by a
  non-trivial group element).
\end{proposition}

The action is generated by the diagonal matrices $$[\exp(2\pi
i\,\overline a_{w_jw_1}),\dots, \exp(2\pi i\,\overline a_{w_jw_n})]$$
for $j=1,\dots,n$, where $\overline a_{ww'}$ is the $(w,w')$-entry of
the inverse matrix
$A(\Gamma)^{-1}$.

In order for this group to act on the complete intersection
singularity, one needs to be able to choose the equations in an
appropriate way.  It turns out that this is sometimes not possible,
even if $\Delta$ satisfies the semigroup condition.  One obtains a new
set of numerical conditions,
which we call the \emph{congruence conditions}; these will guarantee
that the group acts freely (off the origin) on the complete
intersection singularity. Since
the action of the discriminant group depends not only on the
splice diagram, but on the resolution graph, the congruence
conditions depend on the resolution graph.

\begin{conjecture} The resolution graph of a $\Q$-Gorenstein
  singularity with $\Q$-homology sphere link satisfies the congruence
  conditions.
  \end{conjecture}

\begin{example1}
  Returning to our example above, the discriminant group is a group of
  order $16$ whose action on $\C^4$ is generated by the four diagonal
  matrices
$$
[1,-1, i, i], [-1,1,i,i], [i,i,1,-1], [i,i,-1,1].
$$
If this action is to preserve the variety given by equations
(\ref{equations}), the exponents $\alpha$ and $\beta$ must be odd.
This is the congruence condition in this case.  This example is one of
a class of examples (quotient-cusps) for which we have verified our
conjectures.
\end{example1}
\begin{examples}
  For the following example the splice diagram does not
  satisfy the semigroup condition:
$$  
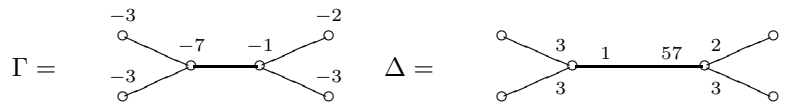
\xymatrix@R=6pt@C=24pt@M=0pt@W=0pt@H=0pt{
\\
&\overtag{\Circ}{-3}{8pt}&&&\overtag{\Circ}{-2}{8pt}\\
{\Gamma=}&&\overtag{\Circ}{-7}{8pt}\lineto[ul]\lineto[dl]\lineto[r]
&\overtag{\Circ}{-1}{8pt}\lineto[ur]\lineto[dr]\lineto[l]&\\
&\overtag{\Circ}{-3}{8pt}&&&\overtag{\Circ}{-3}{8pt}}
\qquad
\xymatrix@R=6pt@C=24pt@M=0pt@W=0pt@H=0pt{
\\
&\Circ&&&&\Circ\\
{\Delta=}&&\Circ\lineto[ul]_(.25){3}\lineto[dl]^(.25){3}\lineto[rr]^(.25){1}^(.75){57}&
&\Circ\lineto[ur]^(.25){2}\lineto[dr]_(.25){3}&\hbox to 0pt{\hglue6pt,\hss}\\
&{\Circ}&&&&{\Circ}\\&~}
$$
since $1$ is not in the semigroup generated by $2$ and $3$.
For the following example the splice diagram satisfies the
semigroup condition, but the resolution graph does not
satisfy the congruence condition:
$$  
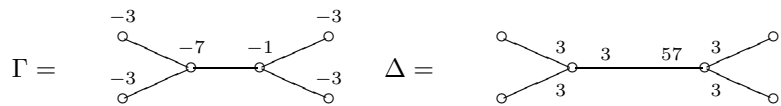
\xymatrix@R=6pt@C=24pt@M=0pt@W=0pt@H=0pt{
\\
&\overtag{\Circ}{-3}{8pt}&&&\overtag{\Circ}{-3}{8pt}\\
{\Gamma=}&&\overtag{\Circ}{-7}{8pt}\lineto[ul]\lineto[dl]\lineto[r]
&\overtag{\Circ}{-1}{8pt}\lineto[ur]\lineto[dr]\lineto[l]&\\
&\overtag{\Circ}{-3}{8pt}&&&\overtag{\Circ}{-3}{8pt}}
\qquad
\xymatrix@R=6pt@C=24pt@M=0pt@W=0pt@H=0pt{
\\
&\Circ&&&&\Circ\\
{\Delta=}&&\Circ\lineto[ul]_(.25){3}\lineto[dl]^(.25){3}\lineto[rr]^(.25){3}^(.75){57}&
&\Circ\lineto[ur]^(.25){3}\lineto[dr]_(.25){3}&\\
&{\Circ}&&&&{\Circ}\\&~\\&~}
$$
\end{examples}

There are two important classes of singularities that are
  topologically determined to be $\Q$-Gorenstein singularities with
  $\Q$-homology sphere links: the rational singularities and (with two
  classes of exceptions) the minimally elliptic singularities of
  Laufer \cite{laufermin}.  In the first non-trivial case we have verified that
  the above purportedly necessary numerical conditions on the graph are
  indeed satisfied.

 \begin{theorem}
     Let $(X,o)$ be a rational or minimally elliptic
 singularity (not simple elliptic or a cusp), whose resolution graph
 has two nodes (vertices of \valence{} at least 3).  Then the splice
 diagram of the resolution graph satisfies the semigroup condition,
 and the congruence conditions are satisfied as well.
 \end{theorem}

When the semigroup and congruence conditions
 are satisfied, one should show that the quotient of the splice
 diagram singularity by the group action has the same topological
 type as $(X,o)$, but we have not yet proved this.

We have a natural notion of equisingular deformation of the above
complete intersections.  We may add to each equation a linear
combination of monomials which for each node $v$ have $v$-weight
greater than or equal to $d_v$, and so that the discriminant group
continues to transform each equation by a character.  These ought
to give equisingular deformations of the original quotient.  In the
case of one node, when dealing with \Brieskorn{} complete intersections,
these are the deformations of weight $\geq 0$.  Our final Conjecture
will say that all  $\Q$-Gorenstein
  singularities with $\Q$-homology sphere link arise in this way.

\begin{conjecture}\label{conjs}
  Let $(X,o)$ be a $\Q$-Gorenstein singularity with $\Q$-homology
  sphere link.  Let $\Delta$ be the splice diagram associated to the
  resolution graph.  Then $\Delta$ satisfies the semigroup condition,
  the resolution graph satisfies the congruence condition, and the
  universal abelian cover of $(X,o)$ is an equisingular deformation of
  the complete intersection of splice diagram form defined above.  The
  geometric genus of the singularity is determined by the topology.
\end{conjecture}

For the last assertion concerning the geometric genus, there are some
deep partial results due to A. N\'{e}methi in \cite{nemethi}.

\section{Results and evidence}

We have verified the conjectures for several examples. The evidence
for our conjectures also includes phenomena that seem to hold in more
general situations than those described by the conjectures.  We
describe now some of this evidence.

\subsection{The topology of the universal abelian cover is determined by
  the splice diagram} Although a resolution tree determines a splice
diagram, the resolution tree cannot be recovered from the splice
diagram. If the splice diagram has just one node (this is the
situation of Theorem \ref{th:weighted homog}) then there are
infinitely many resolution trees sharing the same splice diagram;
otherwise just finitely many. If our conjectures are correct then the
homeomorphism type of the universal abelian cover of the link of a
$\Q$-Gorenstein singularity with $\Q$-homology sphere link should only
depend on the splice diagram of the resolution tree.  This is indeed
true, even without the $\Q$-Gorenstein assumption:
\begin{theorem}
  If two singularities whose links are $\Q$-homology spheres have
  the same splice diagram then the universal abelian covers of their
  links are diffeomorphic.
\end{theorem}

\subsection{Quotient-cusps}
We have succeeded in proving our conjectures in a number of cases.
The quasihomogeneous case was already known by \cite{neumann83}.
A new example is treated in \cite{neumann-wahl00}.  Among the
log-canonical surface singularities, the only non-weighted
homogeneous $\Q$-homology sphere link singularities are the
\emph{quotient-cusps}.  These are $\Z/2$ quotients of certain cusp
singularities; so it is obvious that the universal abelian cover
is also a cusp.  Given a resolution graph of a cusp, it is easy to
tell if it is a complete intersection, by an old result of Karras.
But it requires some extensive argument to show that in our
situation this always occurs.  To actually write down the
equations of the universal abelian cover, and explicit action of
the discriminant group, was originally an extremely difficult
task.  But the discovery of a simpler way to do it led to the
description of the generalized \Brieskorn{} equations plus diagonal
group action considered above.

\subsection{Integral homology sphere links}
In the special case that the link $M$ of our singularity is an
integral homology sphere, that is $H_1(M;\Z)=\{0\}$, our conjectures
become:
\begin{conjecture}
Let $(X,o)$ be a Gorenstein normal surface singularity whose link
is an integral homology sphere.  Then $X$ is a complete
intersection singularity, its splice diagram $\Delta$ satisfies
the semigroup condition, and its equations are equisingular
deformations of those in Section 2.  The Casson invariant of the link
is one-eighth the signature of the Milnor fibre of $X$ (cf.
\cite{neumann-wahl90}).
\end{conjecture}

For a splice diagram with two nodes that represents an integral
homology sphere, as well as some other cases, we can show that our
equations define a singularity whose link has the correct topological
type, and for which the Casson invariant is as conjectured
\cite{neumann-wahl01}.

\end{document}